\documentclass{amsart}

\usepackage{amssymb}
\usepackage{amsmath}
\usepackage{amsthm}
\usepackage{amsbsy}
\usepackage{bm}

\theoremstyle{plain}
\newtheorem{theorem}{Theorem}[section]
\newtheorem{lemma}[theorem]{Lemma}

\newtheorem{proposition}[theorem]{Proposition}

\theoremstyle{definition}
\newtheorem{definition}[theorem]{Definition}

\theoremstyle{remark}

\newcommand{\Z}{\mathbb{Z}}
\newcommand{\R}{\mathbb{R}}
\newcommand{\C}{\mathbb{C}}
\newcommand{\N}{\mathbb{N}}

\newcommand{\del}{\partial}
\newcommand{\s}{\mathfrak{s}}

\begin{document} 

\title{Open manifolds, Ozsvath-Szabo invariants and Exotic $\R^4$'s}

\author{Siddhartha Gadgil}

\address{	Stat Math Unit,\\
		Indian Statistical Institute,\\
		Bangalore 560059, India}

\email{gadgil@isibang.ac.in}

\date{\today}

\subjclass{Primary 57R58; Secondary 53D35,57M27,57N10}

\begin{abstract}
We construct an invariant of certain open four-manifolds using the
Heegaard Floer theory of Ozsvath and Szabo. We show that there is a
manifold $X$ homeomorphic to $\R^4$ for which the invariant is
non-trivial, showing that $X$ is an exotic $\R^4$. This is the first
invariant that detects exotic $\R^4$'s.
\end{abstract}

\maketitle

\section{Introduction}

In this paper, we construct invariants of certain open $4$-manifolds
using the Heegaard Floer theory of Ozsvath and Szabo, and show that our
invariants can detect exotic $\R^4$s. Previous constructions of exotic
$\R^4$'s used indirect arguments to establish exoticity.

Given an $(n+1)$-dimensional field theory, a direct limit construction
can be used to construct an invariant of open $(n+1)$-dimensional
manifolds (which we see in detail later). The subtlety in the case of
Ozsvath-Szabo invariants is that they do not give a field theory, but
satisfy a more complicated composition law. However if we restrict to
a class of cobordisms, which we call \emph{admissible cobordisms}, we
do get a field theory. Using this, we construct our invariants.

Recall that the Ozsvath-Szabo invariants of a smooth, oriented
$3$-manifold $M$ associate homology groups to $M$ equipped with a
$Spin^c$ structure $t$. Further, given a smooth cobordism $W$ between
$3$-manifolds $M_1$ and $M_2$ and a $Spin^c$ structure $\s$ on $W$, we
get an induced map on the groups associated to the restrictions of $\s$
to $M_1$ and $M_2$. To make this into a field theory, one needs a
composition rule for a cobordism $W_1$ from $M_1$ to $M_2$ equipped
with a $Spin^c$ structure $\s_1$ and a cobordism $W_2$ from $M_2$ to
$M_3$ equipped with a $Spin^c$ structure $\s_2$ with
$\s_1|_{M_2}=\s_2|_{M_2}$. However, such $Spin^c$ structures $\s_1$ and
$\s_2$ do not in general uniquely determine a $Spin^c$ structure on the
composition $W=W_1\coprod_{M_2} W_2$ of $W_1$ and $W_2$. We do have a
weaker composition law, where we sum over $Spin^c$ structures on $W$
restricting to $\s_1$ and $\s_2$.

We now find sufficient conditions under which $\s_1$ and $\s_2$
uniquely determine a $Spin^c$ structure $\s$ on $W$. The $Spin^c$
structures on a manifold $X$ are a torseur of $H^2(X,\Z)$. Consider
the Mayer-Vietoris sequence for $W=W_1\cup W_2$
$$\to H^1(W_1)\oplus H^1(W_2)\to H^1(M_2)\overset{\delta}\to H^2(W)\to
H^2(W_1)\oplus H^2(W_2)\to H^2(M_2)$$

From this sequence, it follows that, given $\s_1$ and $\s_2$ as above,
there is a unique $Spin^c$ structure $\s$ on $W$ which restricts to
$\s_1$ and $\s_2$ if and only if the coboundary map $\delta:H^1(M_2)\to
H^2(W)$ is trivial. This is equivalent to the map induced by
inclusions $H^1(W_1)\oplus H^1(W_2)\to H^1(M_2)$ being
surjective. Motivated by this, we make the following definition.

\begin{definition}
A smooth $4$-dimensional cobordism $W$ from $M_1$ to $M_2$ is admissible
if the map induced by inclusion $H^1(W)\to H^1(M_2)$ is surjective.
\end{definition}

We shall see basic properties of such cobordisms in
Section~\ref{cnvx}. We now turn to the corresponding notions for open
manifolds. Let $X$ be an open $4$-manifold which we assume for
simplicity has one end. Let $K_1\subset K_2\subset \dots$ be an
exhaustion of $X$ by compact manifolds and let $M_i=\del K_i$. We
assume here and henceforth (for all exhaustions) that $K_i\subset
int(K_{i+1})$. For $i<j$, let $W_{ij}=K_j-int(K_i)$ be cobordisms from
$M_i$ to $M_j$.

\begin{definition}
The exhaustion $\{K_i\}$ of $X$ is said to be admissible if each
cobordism $W_{ij}$, $i,j\in\N$, $i<j$, is admissible. The manifold $X$
is said to be admissible if it has an admissible exhaustion.
\end{definition}

We shall need to consider the appropriate notion of $Spin^c$ structures
for the ends of $4$-manifolds.

\begin{definition}
An asymptotic $Spin^c$ structure $\s$ on $X$ is a $Spin^c$ structure on
$X-K$ for a compact subset $K\subset X$. Two asymptotic $Spin^c$
structures $\s_1$ and $\s_2$, defined on $X-K_1$ and $X-K_2$, are said
to be equal if there is a compact set $K_0\supset K_1,K_2$ with
$\s_1|_{M-K_0}=\s_2|_{M-K_0}$.
\end{definition}

Given an admissible open $4$-manifold $X$ and an asymptotic $Spin^c$
structure $\s$, we can define invariants of $X$, which we call the
\textit{End Floer Homology}, using direct limits. We shall see in
Section~\ref{inv} that an admissible exhaustion gives a directed
system.

\begin{theorem}
There is an invariant $HE(X,\s)$ which is the direct limit of the
reduced Heegaard Floer homology groups $HF^+_{red}(M_i,\s|_{M_i})$ under
morphisms induced by the cobordisms $W_{ij}$. Furthermore this is
independent of the admissible exhaustion of $X$.
\end{theorem}

We shall also need a \emph{twisted} version of these invariants. Let
$K\subset X$ be a compact set, $\s$ a $Spin^c$-structure on $X-K$ and
$\omega$ a $2$-form on $X-K$. Then we consider the reduced Floer
theory with $\omega$-twisted coefficients (as in~\cite{OZ4}). Once
more we get a directed system whose limit gives an invariant
$\underline{HE}(X,\s)$.

By taking an exhaustion of $\R^4$ by balls, we have the following
proposition.

\begin{proposition}
For the unique asymptotic $Spin^c$ structure $\s$ on $\R^4$ (and any $2$-form $\omega$ on $\R^4-K$ with $K$ compact), we have
$\underline{HE}(\R^4,\s)=0$.
\end{proposition}

Our main result is that there are manifolds homeomorphic to $\R^4$ but
with non-vanishing end Floer homology.

\begin{theorem}\label{exot}
There is a $4$-manifold $X$ homeomorphic to $\R^4$ such that there is
a compact set $K\subset X$, a $Spin^{\C}$ structure $\s$ on $X-K$ and a
closed $2$-form $\omega$ on $X-K$ with $\underline{HE}(X,\s)\neq 0$
with $\omega$-twisted coefficients.
\end{theorem}

Thus, $X$ is an exotic $\R^4$. Previous constructions of exotic
$\R^4$'s used indirect arguments to show that they are exotic. The
\textit{End Floer homology} is the first invariant that detects exotic
$\R^4$'s.

\section{Admissible cobordisms and admissible ends}\label{cnvx}

We henceforth assume that all our manifolds are smooth and oriented
and all cobordisms are compact and $4$-dimensional. By $W:M_1\to M_2$
we mean a smooth cobordism from the closed $3$-manifold $M_1$ to the
closed $3$-manifold $M_2$. Given $W_1:M_1\to M_2$ and $W_2:M_2\to
M_3$, $W_2\circ W_1$ denotes the composition of the cobordisms $W_1$
and $W_2$.

In this section we prove some simple
results concerning admissible cobordisms and admissible ends.

\begin{lemma}\label{comp}
Suppose $W_1:M_1\to M_2$ and $W_2:M_2\to M_3$ are admissible
cobordisms, then $W=W_2\circ W_1$ is admissible.
\end{lemma}
\begin{proof}
We need to show that the map $H^1(W)\to H^1(M_3)$ induced by inclusion
is surjective. This is the composition of maps $H^1(W)\to H^1(W_2)$
and $H^1(W_2)\to H^1(M_3)$ induced by inclusion, with the latter
surjective by hypothesis. We shall show that the map $H^1(W)\to
H^1(W_2)$ is surjective.

Let $\alpha\in H^1(W_2)$ be a class. Let $i_j:M_2\to W_j$, $j=1,2$, be
inclusion maps. Consider the Mayer-Vietoris sequence
$$\dots\to H^1(W)\to H^1(W_1)\oplus H^1(W_2)\overset{i_1^*+i_2^*}\to
H^1(M_2)\to\dots $$ 

By admissibility of $W_1$, there is a class $\beta\in H^1(W_1)$ with
$i_1^*(\beta)=i_2^*(\alpha)$. Hence the image of the class
$(-\beta,\alpha)\in H^1(W_1)\oplus H^1(W_2)$ in $H^1(M_2)$ is zero,
and so $(-\beta,\alpha)$ is the image of a class $\varphi\in
H^1(W)$. In particular $\alpha$ is the image of $\varphi$ under the
map induced by inclusion.
\end{proof}

\begin{lemma}\label{comp2}
Suppose $W_1:M_1\to M_2$ and $W_2:M_2\to M_3$ are cobordisms with
$W=W_2\circ W_1$ admissible. Then $W_2$ is admissible.
\end{lemma}
\begin{proof}
By hypothesis the map $H^1(W)\to H^1(M_3)$ is surjective. This factors
through the map $H^1(W_2)\to H^1(M_3)$, which must also be surjective.
\end{proof}

We need criteria for when cobordisms corresponding to attaching
handles are admissible.

\begin{lemma}\label{handl} 
Let $M=M_1$ be a $3$-manifold, $W$ the cobordism corresponding to a
handle addition and $M_2$ the other boundary components of $W$. The
following hold.
\begin{enumerate}
\item A product cobordism is admissible.
\item The cobordism corresponding to attaching a $1$-handle to a
closed $3$-manifold $M$ is admissible.
\item If $K$ is a knot in a closed $3$-manifold which represents a
primitive, non-torsion element in $H_1(M)$, then the cobordism
corresponding to attaching a $2$-handle to $M$ along $K$ is admissible.
\end{enumerate}
\end{lemma}
\begin{proof}
We shall show that the map induced by the inclusion from $H_1(M_2)$ to
$H_1(W)$ is an isomorphism in each case. As the map on cohomology is
the adjoint of this map, it follows that it is a surjection.

The case of a product cobordism is immediate. In the second case we
see that $H_1(M_2)=H_1(W)=H_1(M)\oplus \Z$ with the isomorphism
induced by inclusion. In the third case we have $H_1(M)=H\oplus\Z$,
with $[K]$ generating the $\Z$ component and $H$ isomorphic to the
homology of the $3$-manifold obtained by surgery about $K\subset M$. It
is easy to see that $H_1(W)=H_1(M_2)=H$.
\end{proof}

Now let $X$ be an open manifold and let $K_1\subset K_2\subset\dots $
be an exhaustion of $X$ and $M_i$ and $W_{ij}$ be as before.

\begin{lemma}
The exhaustion $\{K_i\}$ is admissible if and only if each of the
manifolds $K_{j+1}-int(K_j)$ is admissible.
\end{lemma}
\begin{proof}
Each $W_{ij}$ is the composition of cobordisms $K_{j+1}-int(K_j)$. The
result follows by Lemmas~\ref{comp} and~\ref{comp2}.
\end{proof}

Thus, if $X$ is obtained from a compact manifold $K$ by attaching
handles as in Lemma~\ref{handl} then $X$ is admissible. Our examples
of exotic $\R^4$s will be of this form.

It is immediate from the definition that for any admissible exhaustion
$K_i$, the exhaustion obtained by passing to a subsequence $K_{i_j}$
is admissible. To show independence of our invariants under
exhaustions, we need the following lemma.

\begin{lemma}\label{refin}
Let $K_1\subset L_1\subset K_2\subset L_2\dots$ be an exhaustion of
$X$ with $K_1\subset K_2\subset \dots$ and $L_1\subset L_2\subset
\dots$ admissible exhaustions. Then the exhaustion $L_1\subset
K_2\subset L_2\subset K_3\dots$ is admissible.
\end{lemma}
\begin{proof}
It suffices to show that the cobordisms $K_{j+1}-int(L_j)$, $j\geq 1$
and $L_j-int(K_j)$, $j\geq 2$ are admissible. This follows from
Lemma~\ref{comp2} as the cobordisms $K_{j+1}-int(K_j)$ and
$L_{j+1}-int(L_j)$ are admissible and we have
$K_{j+1}-int(K_j)=(K_{j+1}-int(L_j))\circ (L_{j}-int(K_j))$ and
$L_{j+1}-int(L_j)=(L_{j+1}-int(K_j))\circ (K_{j}-int(L_j))$.
\end{proof}

\section{Invariants for admissible ends}\label{inv}

We are now ready to define our invariants for an admissible open
$4$-manifold $X$. We shall construct invariants based on reduced
Heegaard Floer theory $HF_{red}^+$. First we recall some facts about
Ozsvath-Szabo theory.

Associated to each closed, oriented $3$-manifold $M$ and $Spin^c$
structure $t$ on $M$ we have abelian groups $HF^+(M,t)$,
$HF^-(M,t)$ and $HF^{\infty}(M,t)$ that fit in an exact sequence
$$\dots\to HF^-(M,t)\to HF^{\infty}(M,t)\to HF^+(M,t)\to\dots$$

Further, a cobordism $W:M_1\to M_2$ with a $Spin^c$ structure $\s$ on
$W$ such that $t_i=s|_{M_i}$ induces homomorphisms $F_{W,\s}$ on these
abelian groups which commute with the maps in the above exact
sequence.

The group $HF_{red}^+(M,t)$ is defined as the quotient of $HF^+(M,t)$
by the image of $HF^{\infty}(M,t)$. This is isomorphic to the
kernel $HF_{red}^-(M,t)$ of the map from $HF^-(M,t)$ to
$HF^{\infty}(M,t)$. Further, given a cobordism $W:M_1\to M_2$ with a
$Spin^c$ structure $\s$ on $W$ such that $t_i=s|_{M_i}$, we get an
induced homomorphism  on the abelian groups
$F_{W,\s}:HF_{red}^+(M_1,t_1)\to HF_{red}^+(M_2,t_2)$ induced by the
corresponding homomorphism on $HF^+$ as the image of
$HF^{\infty}(M_1,t)$ is contained in $HF^{\infty}(M_2,t)$. This
homomorphism is well defined up to choice of sign. We shall denote the
above cobordism with its $Spin^c$ structure by $(W,\s):(M_1,t_1)\to
(M_2,t_2)$.

Further, if $(W_1,\s_1):(M_1,t_1)\to (M_2,t_2)$ and
$(W_2,\s_2):(M_2,t_2)\to (M_3,t_3)$, with $W=W_2\circ W_1$, we have the
composition formula
$$F_{W_2,\s_2}\circ F_{W_1, s_1}=\sum _{s|_{W_i}=s_i} \pm F_{W,\s}$$

We shall consider the special case when $W_1$ is admissible.

\begin{lemma}\label{fact}
If $W_1$ is admissible then there is a unique $Spin^c$ structure $\s$ on
$W$ with $\s|_{W_i}=s_i$. For this $Spin^c$ structure $F_{W_2,\s_2}\circ
F_{W_1, s_1}= \pm F_{W,\s}$
\end{lemma}
\begin{proof}
Recall that $Spin^c$ structures are a torseur of
$H^2(\cdot,\Z)$. Consider the Mayer-Vietoris sequence for $W=W_1\cup
W_2$
$$\to H^1(W_1)\oplus H^1(W_2)\to H^1(M_2)\overset{\delta}\to H^2(W)\to
H^2(W_1)\oplus H^2(W_2)\to H^2(M_2)$$

By admissibility the map $H^1(W_1)\oplus H^1(W_2)\to H^1(M_2)$ is a
surjection, hence $H^2(W)\to H^2(W_1)\oplus H^2(W_2)$ is an
injection. This shows uniqueness of the $Spin^c$ structure. As
$\s_1|_{M_2}=t_2=s_2|_{M_2}$, existence follows from the same exact
sequence.

The second statement follows from the first using the composition formula.
\end{proof}

For an admissible exhaustion, it follows that we get a directed system
of abelian groups up to sign.  We next see that we can choose signs to
get a directed system, and the direct limit of the system does not
depend on the choice of signs.

\begin{lemma}\label{sign}
Assume $A_i$ is a sequence of Abelian groups and maps $f_{ij}:A_i\to
A_j$, such that for $i<j<k$, $f_{ik}=\pm f_{jk}\circ f_{ij}$. Then we
can choose $g_{ij}=\pm f_{ij}$ such that we get a directed
system. Furthermore the limit is independent, up to isomorphism, of
the choices.
\end{lemma}
\begin{proof}
Let $g_{1j}=f_{1j}$. For $i<j$, the composition law
$g_{1j}=g_{ij}\circ g_{1i}$ uniquely determines sign of $g_{ij}=\pm
f_{ij}$, and such a $g_{ij}$ exists as $f_{1j}=\pm f_{ij}\circ
f_{1i}$. It is easy to see that this gives a directed system.

For a different choice the maps $g_{1j}$ are replaced by
$g'_{ij}=\epsilon_j g_{1j}$, $\epsilon_j=\pm 1$. We get in general a
different directed system, with the groups $A_i$ . However, using the
isomorphisms $\epsilon_i:A_i\to A_i$ (i.e., $x\mapsto \epsilon_i\times
x$ for $x\in A_i$), we get an isomorphism of directed systems. Hence
the limits are isomorphic.
\end{proof}

\begin{definition}
The End Floer homology $HE(X,\s)$ is the direct limit of the directed
system constructed above.
\end{definition}

\begin{proposition}\label{welldef}
The End Floer homology is independent of the admissible exhaustion chosen. 
\end{proposition}
\begin{proof}
By elementary properties of direct limits, the limit does not change
on passing to a subsequence of an exhaustion. Given two admissible
exhaustions $K_1\subset K_2\subset \dots$ and $L_1\subset L_2\subset
\dots$, by passing to subsequences we can assume that $K_1\subset
L_1\subset K_2\subset L_2\subset \dots$ for the two exhaustions. By
Lemma~\ref{refin} the exhaustion $L_1\subset K_2\subset L_2\subset
K_3\dots$ is admissible. As $L_1\subset L_2\subset \dots$ and
$K_2\subset K_3\subset \dots$ are subsequences of this exhaustion, the
direct limits for the exhaustions $K_1\subset K_2\subset \dots$ and
$L_1\subset L_2\subset \dots$ are the same (as they are both
isomorphic to the direct limit corresponding to the exhaustion
$L_1\subset K_2\subset L_2\subset K_3\dots$).
\end{proof}

We see that this depends only on the diffeomorphism class of the end of $X$. More precisely, we have the following.

\begin{proposition}\label{enddiff}
Suppose $X$ and $Y$ are admissible smooth $4$-manifolds  and $K\subset X$ and $L\subset Y$ are compact sets so that there is a diffeomorphism $f:X-K\to Y-L$. Then the End Floer homology groups of $X$ and $Y$ are isomorphic.
\end{proposition}
\begin{proof}
Consider an admissible
exhaustion $K_1\subset K_2\subset \dots$ with $K\subset K_1$. We define an exhaustion $L_1\subset L_2\subset \dots$ of $Y$ by $L_i=L\cup f(K-K_i)$. The map $f$ induces isomorphisms between the terms of the directed systems corresponding to the two exhaustions. Thus, the End Floer homology groups, which are the limits of these directed systems, are isomorphic.
\end{proof}

We consider the $\omega$-twisted version of this as in~\cite{OZ4}. Let
$K\subset X$ be a compact manifold and $\omega$ a $2$-form on
$X-K$. We call such a $2$-form $\omega$ on $X-K$, for $K$ compact, an asymptotic $2$-form. Given two closed $2$-forms $\omega_i$, $i=1,2$, on the complements $X-K_i$ of smooth compact sets $K_i$, $1=1,2$, we say that $\omega_1$ and $\omega_2$ are asymptotically cohomologous if, for some compact set $K$, $K_i\subset K$ for $i=1,2$, the restrictions of the forms are cohomologous on $X-K$. We can thus speak of asymptotic cohomology classes of asymptotic $2$-forms.

We consider an admissible exhaustion with the first term $K_1$
satisfying $K\subset K_1$. For this, we can define the twisted groups
$\underline{HF}_{red}^+(M_i,t_i)$ and homomorphisms associated to
$W_{ij}$ which are well defined up to sign and multiplication by
powers of $T$. For any composition $W=W_2\circ W_1$ associated with
the exhaustion as above, the coboundary map $\delta:H^1(M_2)\to
H^2(W)$ is zero. It follows by the composition rule for
$\omega$-twisted coefficients that we have a directed system up to
multiplication by powers of $T$ and sign. As in Lemma~\ref{sign}, we
can make choices for the homomorphisms to get a directed system and
the direct limit is independent of the choices.

The direct limit is the End Floer homology $\underline{HE}(X,\s)$ with
$\omega$-twisted coefficients. The following propositions are ananlogous to Propositions~\ref{welldef} and~\ref{enddiff}.

\begin{proposition}
For an aymptotic $2$-form $\omega$, the $\omega$-twisted End Floer homology is independent of the choice of admissible exhaustion. 
\end{proposition}

\begin{proposition}
Let $X$ and $Y$ are smooth $4$-manifolds with admissible ends and $\omega_X$ and $\omega_Y$ are asymptotic $2$-forms on $X$ and $Y$. If there are compact sets $K\subset X$ and $L\subset Y$, with $\omega_X$ and $\omega_Y$ defined on $X-K$ and $Y-L$, and a diffeomorphism $f:X-K\to Y-L$ so that $\omega_X$ is asymptotically cohomologous to $f^*(\omega_Y)$, then the End Floer homology with $\omega_X$-twisted coefficients of $X$ is homologous to the End Floer homology with $\omega_Y$ twisted coefficients of $Y$.
\end{proposition}

\section{Exotic $\R^4$'s}

We now construct a manifold $X$ homeomorphic to $\R^4$ with
$\underline{HE}(X)\neq 0$. This is done by first constructing a convex
symplectic manifold $W$ with one convex boundary component $N_0$ and
one convex end and then gluing a compact manifold $Y$ to $W$ along
$N_0$.

\subsection{Construction of $X$}
Let $K$ be a non-trivial slice knot in $S^3$ and let $N$ be obtained
by $0$-frame surgery about $K$. Then $N$ admits a taut
foliation by~\cite{Ga}, and hence $N\times [0,1]$ admits a symplectic structure with both
ends convex by~\cite{ET}. The symplectic structure induces a contact structure $\xi$ on $N$. We shall construct a symplectic manifold $Q$ with one concave boundary component contactomorphic to $(N,\xi)$ and one convex end. The manifold $W$ is obtained by gluing $Q$ to $N\times [0,1]$.

Let $P$ be the manifold obtained by attaching a $2$-handle $H$ to $N\times
\{1\}$ corresponding to the surgery cancelling the $0$-frame surgery
about $K$. The manifold $P$ has boundary $S\cup N_0$ with
$N_0=N\times \{0\}$ and $S$ a $3$-sphere. Let $P_0$ be $P-S$. Then $P_0$ has one boundary component, which is diffeomorphic to $N$, and one end.

\begin{lemma}\label{casson}
There is a symplectic manifold $Q$ properly homotopy equivalent to $P_0$ so that the end of $Q$ is convex and the boundary component identified with $N_0$ is concave with induced contact structure $\xi$. 
\end{lemma}
\begin{proof}
We construct $Q$ as a Stein cobordism as in~\cite{EH}. Firstly, by a theorem of Eliashberg~\cite{El1} (Lemma~2.2 in~\cite{EH}), there is a Stein cobordism from $(N,\xi)$ to itself, which is thus a Stein structure on $N\times [0,1]$ with $N\times\{0\}$ a concave boundary component and $N\times \{1\}$ a convex boundary component. We construct the manifold $Q$  by attaching $1$-handles and $2$-handles starting with the convex boundary component, with the $2$-handles attached with framing $1$ less than the Thurston-Bennequin framing (we call this Legendrian handle addition). By Eliashberg's characterisation of Stein domains~\cite{El2} (see also~\cite{El3} and~\cite{Go}), $Q$ is Stein.

The $1$-handles and $2$-handles are attached as in Theorem~3.1 of~\cite{Go}, so that the handle $H$ is replaced by a Stein Casson handle. Specifically, by taking a Legendrian representative of $\kappa=\del H$, we can perform Legendrian handle addition about $\kappa$ but with incorrect framing, differing from that of $H$ by an integer $k$. If we attach a handle to $\kappa$ with this farming but with $k$ self-plumbings (a so called \emph{kinky handle}), then the self-intersection pairing coincides with that obtained by attaching $H$. As in~\cite{Go} (where there is an explicit construction in Figure~22), one can attach $1$-handles and Legendrian $2$-handles to obtain a Stein manifold diffeomorphic to that obtained by attaching a $2$-handle with $k$   self-plumbings to $\kappa$ so that we have the same intersection pairing as adding $H$.

Thus, we obtain a Stein cobordism with the same intersection pairing as attaching the handle $H$, but with non-trivial fundamental group. By a lemma of Casson, we can find a family of curves on the boundary of the attached kinky handle, hence the convex boundary of the Stein cobordism, so that attaching $2$-handles to these curves (with appropriate framing) gives the manifold obtained on attaching $H$. As before, we can instead attach kinky handles to obtain a Stein cobordism. 

Iterating this procedure gives a non-compact Stein cobordism $Q$ with one concave boundary component and one convex end, which is diffeomorphic to the manifold obtained by attaching a Casson handle in place of $H$. As Casson handles are properly homotopy equivalent to the interiors of handles, $Q$ is properly homotopy equivalent to $P_0$.

\end{proof}

Let $W$ be the symplectic manifold obtained by gluing $N\times [0,1]$ with its sympectic structure obtained by the Gabai-Eliashberg-Thurston theorem, to the symplectic manifold $Q$, with $N\times \{1\}$ identified with the (concave) boundary of $Q$. Observe that $W$ is simply-connected as the Casson handle corresponding to the $2$-handle $H$ is attached along the
meridian of $K$, which normally generates $\pi_1(N)$. Also observe
that in the proof of Lemma~\ref{casson}, following Theorem~3.1 of~\cite{Go}, the handles attached
are as in Lemma~\ref{handl}, and hence the corresponding exhaustion is
admissible.

Next, let $Y'$ be obtained from $B^4$ by attaching a $2$-handle along
$K$ with framing~$0$. Then $\del Y'=N$. As $K$ is slice, the generator
of $H_2(Y)=\Z$ can be represented by an embedded sphere $\Sigma$. Let
$Y$ be obtained from $Y'$ by performing surgery along $\Sigma$. Glue
$W$ to $Y$ along $\del Y=N=N\times\{0\}$ to obtain $X$.

By a Mayer-Vietoris argument, $X$ has the homology of $\R^4$. Further,
as $\pi_1(Y)$ is normally generated by a meridian of $K$, to which a
Casson handle is attached, $\pi_1(X)=1$. Finally, the end of $X$ is
properly homotopic to the end of $P_0=P-S$, and hence $Y$ is simply-connected at
infinity.  Thus $Y$ is homeomorphic to $\R^4$ by Freedman's theorem~\cite{Fr}.

\subsection{Non-Vanishing of End Floer homology}
Finally, we show that the End Floer homology for $X$ does not
vanish. Consider the exhaustion of $X$ with $K_1=Y$, hence $M_1=N$ and
$K_2$, $K_3$, \dots being the level sets after attaching successive
handles as above. Note that $X-K_1$ is symplectic with symplectic form
$\omega$, and each of the cobordisms $W_{1j}$ is a convex symplectic
manifold with two convex boundary components $M_1$ and $M_j$. Hence $W_{1j}$
embeds in a symplectic $4$-manifold $Z=X_1\cup W_{1j}\cup X_j$ with
both components of $Z-W_{1j}$ having $b_2^+>0$ by results of
Eliashberg~\cite{El} and Kronheimer-Mrowka~\cite{KM}. Here $X_1$ and
$X_j$ are manifolds with boundaries $M_1$ and $M_j$, respectively.

We shall consider $\omega$-twisted coefficients and the $Spin^c$
structure $\s$ associated to $\omega$. Recall that $\omega$-twisted
coefficients are coefficients determined by $\omega$ as follows: for a
$3$-manifold $P\subset M$, we consider $\Z[\R]$ as a module over
$\Z[H^1(N,\Z)]$ via the ring homomorphism $[\gamma]\mapsto T^{\int_N
[\gamma]\wedge\omega}$. Ozsvath and Szabo show that we have induced
maps with $\omega$-twisted coefficients satisfying an appropriate
composition formula. By an application of Stokes theorem, we deduce
the relation
$$\int_N [\gamma]\wedge\omega=\int_Z \delta[\gamma]\wedge\omega$$

Let $t_i$ be the $Spin^c$ structure on $M_i$ induced by $\s$.  We first
construct an element $x_1\in \underline{HF}^+(M_1,t_1)$ whose image
$z_1\in \underline{HF}_{red}^+(M_1,t_1)$ will be shown to have
non-zero image in the direct limit giving the End Floer homology.

Let $P\subset X_1$ be an admissible cut in the terminology of Ozsvath
and Szabo. Then as $\delta H^1(P)=0$, $\omega$-twisted coefficients
coincide with untwisted coefficients(as $\int_P
[\gamma]\wedge\omega=\int_Z \delta[\gamma]\wedge\omega=0$). Let the
closures of the components of $X_1-P$ be $U$ and $V$, with
$M_1\subset\del V$. Let $B_1\subset U$ be a ball. As in the
construction of the closed $4$-manifold invariants, we obtain an
element $\xi\in \underline{HF}^+(P,\s)={HF}^+(P,\s)$ as the image of the
generator of $HF^-(S^3)$ using the isomorphism between $HF^-_{red}$
and $HF^+_{red}$. We define $x_1$ to be the image
$\underline{F}_V(\xi)$ of $\xi$ in $\underline{HF}^+(M_1,t_1)$ under
the map induced by the cobordism $V$ and let $z_1$ be its image in
reduced Floer homology.

Let $x_j\in\underline{HF}^+(M_j,t_j)$ be the image of $x_1$ under the
cobordism induced by $W_{1j}$ and let $z_j\subset
\underline{HF}_{red}^+(M_j,t_j)$ be corresponding image of $z_1$. 

\begin{lemma}
For every $j\geq 0$, $z_j\neq 0$.
\end{lemma}
\begin{proof}
Let $j>1$ be fixed.  Let $W=W_{1j}\cup X_j$ and let $B_2$ be a ball in
$X_j$. We shall show that the image of $x_1$ in $HF^+(S^3,\s_0)$ under
the map induced by $W-B_2$ is non-zero. 

\begin{lemma}
The image $\underline{F}_{W-B_2}(x_1)$ of $x_1$ in
$HF^+(S^3,\s_0)$ under the map induced by $W-B_2$ is non-zero.
\end{lemma}
\begin{proof}
Our proof is based on the proof of Theorem~4.2 in~\cite{OZ4}. We use the product formula with $\omega$-twisted coefficients
$$\sum_{\eta\in H^1(M_1,\Z)} \Phi_{M,\s+\delta\eta}T^{<\omega\cup
c_1(s+\delta\eta),[M]>}=\underline{F}_{W-B_2}\circ
\underline{F}_V(\xi)=\underline{F}_{W-B_2}(x_1)$$

Thus it suffices to show that the left hand side does not vanish. By
results of Ozsvath and Szabo on the closed four-manifold invariants for symplectic manifolds (as in~\cite{OZ4}, Theorem~4.2), the lowest order term of the left hand side, which is a polynomial in $T$, is $1$. It follows that
$\underline{F}_{W-B_2}(x_1)\neq 0$, completing the proof.
\end{proof}

Now, by Lemma~\ref{fact}, as $W_{1j}$ is admissible, this factors
through the map induced by $W_{1j}$, and hence the image of $x_j$ in
$HF^+(S^3,\s_0)$ is non-zero. But as the cobordism $X_j-int(B_2)$ has
$b_2^+>0$, the induced map on $\underline{HF}^{\infty}$ is zero. It
follows that $x_j$ is not in the image of
$\underline{HF}^{\infty}(M_i,t_i)$, i.e. $z_j\neq 0$, as claimed.
\end{proof}

Thus, the End Floer homology of $X$ does not vanish. We have seen that
$X$ is homeomorphic to $\R^4$. This completes the proof of
Theorem~\ref{exot}.\qed

 \bibliographystyle{amsplain}

\end{document}